\documentclass[12pt]{article}%
\usepackage{amsmath}
\usepackage{amsfonts}
\usepackage{amssymb, amsfonts, amsmath, nicefrac}
\usepackage{color}
\usepackage{amssymb}
\usepackage{graphicx}%
\providecommand{\U}[1]{\protect\rule{.1in}{.1in}}
\newcommand{\cg}{{\mathfrak g}}

\newtheorem{theorem}{Theorem}[section]

\newtheorem{proposition}{Proposition}[section]
\newtheorem{example}{Example}[section]

\newtheorem{assumption}{Assumption}[section]

\newcommand{\F}{{\cal F}}

\newcommand{\be}{\begin{equation}}
\newcommand{\ee}{\end{equation}}

\def\ess {{\rm ess}\sup}

\def\dst {{\rightsquigarrow}}

\def\bbe{{\Bbb{E}}}

\begin{document}

\title{Statistical Limit Theorems in Distributionally Robust Optimization}
\author{
}
\maketitle

\begin{center}%
\begin{tabular}
[c]{ccc}%
\begin{tabular}
[c]{c}%
Jose Blanchet\\
Department of Management Science\\
and Engineering\\
Stanford University\\
\texttt{jose.blanchet@stanford.edu}%
\end{tabular}
&  &
\begin{tabular}
[c]{c}%
Alexander Shapiro\\
Georgia Institute of Technology\\
Atlanta, Georgia 30332-0205, USA,\\
\texttt{ashapiro@isye.gatech.edu}\\
\end{tabular}
\end{tabular}

\end{center}

\begin{abstract}
The goal of this paper is to develop methodology for the systematic analysis
of asymptotic statistical properties of data driven DRO formulations based on
their corresponding non-DRO counterparts. We illustrate our approach in
various settings, including both phi-divergence and Wasserstein uncertainty
sets. Different types of asymptotic behaviors are obtained depending on the
rate at which the uncertainty radius decreases to zero as a function of the
sample size and the geometry of the uncertainty sets.

\end{abstract}

\section{Introduction}

The statistical analysis of Empirical Risk Minimization (ERM) estimators is a
well investigated topic both in statistics (e.g., \cite{vandervaart}) and
stochastic optimization (e.g., \cite{SDR}). In recent years, there has been
significant interest in the investigation of distributionally robust
optimization (DRO) estimators (e.g., \cite{RM19}). The goal of this paper is to develop
methodology for the study of asymptotic statistical properties of data driven
DRO formulations based on their corresponding non-DRO counterparts.

Our objective is to illustrate the main conceptual strategies for the
statistical development, emphasizing qualitative features, for instance, the
different types of behavior  arising from the interaction between the
distributional uncertainty size and the sample size, while keeping the
discussion easily accessible. Consequently, in order to keep the discussion
easily accessible, we do not necessarily focus on the most general assumptions
to apply our results.

To set the stage, let us introduce some notation. We use ${\mathfrak{P}%
}(\mathcal{S})$ to denote the set of Borel probability measures supported on a
closed (nonempty) set $\mathcal{S}\subset{\mathbb{R}}^{d}$. Let $X_{1}%
,...,X_{n}$ be a sequence of independent identically distributed (i.i.d.)
random vectors viewed as realizations (or i.i.d. copies) of random vector $X$
having distribution $P_{\ast}\in{\mathfrak{P}}(\mathcal{S})$. Consider the
corresponding empirical measure $P_{n}=n^{-1}\sum_{i=1}^{n}\delta_{X_{i}}$,
where $\delta_{x}$ denotes the Dirac measure of mass one at the point
$x\in{\mathbb{R}}^{d}$. The sample mean of a function $\psi:\mathcal{S}%
\rightarrow{\mathbb{R}}$ is ${\mathbb{E}}_{P_{n}}[\psi(X)]=n^{-1}\sum
_{i=1}^{n}\psi(X_{i})$. By the Strong Law of Large Numbers we have that
${\mathbb{E}}_{P_{n}}[\psi(X)]$ converges with probability one (w.p.1) to
${\mathbb{E}}_{P_{\ast}}[\psi(X)]$, provided the expectation ${\mathbb{E}%
}_{P_{\ast}}[\psi(X)]$ is well defined\footnote{Throughout our discussion, every
function whose expectation is considered will be assumed to be Borel
measurable, so we will not be concerned with making this assumption repeatedly.}.

By the Central Limit Theorem,
\[
n^{1/2}({\mathbb{E}}_{P_{n}}\left[  \psi\left(  X\right)  \right]
-{\mathbb{E}}_{P_{\ast}}[\psi\left(  X\right)  ])\dst N\left(  0,\sigma
^{2}\right)  ,
\]
where \textquotedblleft$\dst$\textquotedblright\ denotes the weak convergence
(converges in distribution) and $N\left(  0,\sigma^{2}\right)  $ represents the
normal distribution  with mean zero and variance $\sigma^{2}=\mathrm{Var}%
_{P_{\ast}}[\psi(X)]$, provided this variance is finite.

We consider a loss function of the form $l:{\mathbb{R}}^{d}\times
\Theta\rightarrow{\mathbb{R}}$, with $\Theta\subset{\mathbb{R}}^{m}$ being the
parameter space. Unless stated otherwise, we assume that the set $\Theta$ is
\emph{compact} and $l(x,\theta)$ is \emph{continuous} on $\mathcal{S}%
\times\Theta$.
We define
\begin{equation}
f_{n}\left(  \theta\right)  :={\mathbb{E}}_{P_{n}}\left[  l\left(
X,\theta\right)  \right]  \text{ \ and }f\left(  \theta\right)  :={\mathbb{E}%
}_{P_{\ast}}\left[  l\left(  X,\theta\right)  \right]  . \label{eq-0}%
\end{equation}
So, the standard ERM formulation takes the form
\begin{equation}
\min_{\theta\in\Theta}f_{n}\left(  \theta\right)  , \label{eq-1}%
\end{equation}
and viewed as an empirical counterpart of the \textquotedblleft
true\textquotedblright\ (or limiting) form
\begin{equation}
\min_{\theta\in\Theta}f\left(  \theta\right)  . \label{eq-2}%
\end{equation}
The statistical properties such as consistency and asymptotic normality of the
ERM estimates have been widely studied in significant generality as the sample
size $n\rightarrow\infty$. These types of results hold under structural
properties of the function $f\left(  \cdot\right)  $ and natural stability
assumptions (to be reviewed) which guarantee a functional Central Limit
Theorems for $f_{n}\left(  \cdot\right)  $. Our goal is to present a
development that is largely parallel for the associated distributionally
robust counterpart to (\ref{eq-1}).

More precisely, (\ref{eq-1}) can be endowed with distributional robustness by
defining a set of probability measures, referred to as the \emph{ambiguity
set}, ${\mathfrak{M}}_{\delta}(P_{n})\subset{\mathfrak{P}}(\mathcal{S})$,
which are seen as \textquotedblleft reasonable\textquotedblright\ (according
to some criterion) perturbations of the empirical measure. The parameter
$\delta\geq0$ is the uncertainty size and the family of sets $\{{\mathfrak{M}%
}_{\delta}(P_{n}):\delta\geq0\}$ is typically nondecreasing in $\delta$ (in
the inclusion partial order sense).
The ambiguity set can be defined around any reference probability measure, but
unless stated otherwise, we will center the ambiguity set around $P_{n}$. In
this paper we deal with ambiguity sets of the form
\begin{equation}
{\mathfrak{M}}_{\delta}(P_{n}):=\{P\in{\mathfrak{P}}(\mathcal{S}%
):D(P,P_{n})\leq\delta\}, \label{ambig-1}%
\end{equation}
where $D(Q,P)$ is a divergence between $Q,P\in{\mathfrak{P}}(\mathcal{S})$.
Specifically, we consider the phi-divergence and Wasserstein distance cases.

In order to state the DRO version of \eqref{eq-1} we define
\begin{equation}
\mathcal{F}_{n}(\theta,\delta_{n}):=\sup_{P\in{\mathfrak{M}}_{\delta_{n}%
}(P_{n})}{\mathbb{E}}_{P}\left[  l\left(  X,\theta\right)  \right]  ,
\label{opt-3}%
\end{equation}
where $\delta_{n}$ is a monotonically decreasing sequence tending to zero as
$n\rightarrow\infty$. The DRO version of \eqref{eq-1} takes the form
\begin{equation}
\min_{\theta\in\Theta}\mathcal{F}_{n}(\theta,\delta_{n}). \label{eq-3}%
\end{equation}
The aim of this paper is to investigate asymptotic statistical properties of
the optimal value and optimal solutions of the DRO problem \eqref{eq-3}. In
particular, under natural assumptions (to be discussed), we will show that
both in phi-divergence and Wasserstein DRO formulations, there are
\textit{typically} (but not always) \textit{three types of cases} involving
the limiting asymptotic statistics depending on the rate of convergence of
$\delta_{n}$ to zero. These can be seen both in terms of the value function
error%
\[
\min_{\theta\in\Theta}\mathcal{F}_{n}(\theta,\delta_{n})-\min_{\theta\in
\Theta}f\left(  \theta\right)  ,
\]
and the optimal solution error (assuming it is unique for the limiting version
of the problem and sufficient regularity conditions are in place).

Intuitively, if $\delta_{n}$ is smaller than a certain (to be characterized)
\textit{critical rate} relative to the canonical parametric statistical error
rate $n^{-1/2}$, then the DRO effect is negligible compared to the statistical
error implicit in a sample of size $n$. If $\delta_{n}$ decreases to zero
right at the critical rate, the DRO effect is comparable with this statistical
error and can be quantified in the form of an asymptotic bias. If $\delta_{n}$
is bigger than the critical rate, the DRO\ effect overwhelms the statistical
noise. These critical rates depend on the sensitivity of the optimal value
function with respect to a small change in the size of uncertainty.

Our objective is to provide accessible principles that can be used to obtain
explicit limiting distributions for the errors, both for value functions and
optimizers, when $\delta_{n}\rightarrow0$ in these \textit{three cases}; see
Theorems \ref{Thm_GP_DRO_OV} and    \ref{Thm_GP_DRO_OS} for general
principles and Theorems \ref{th-asymdro} and \ref{th_value_DRO_W} for the
application to these principles to value functions of phi-divergence and
Wasserstein DRO, respectively; and Theorems \ref{thm_DRO_OS_phi} and  \ref{Thm_DRO_OS_W} for the corresponding application to phi-divergence and
Wasserstein DRO optimal solutions, respectively.

It is important to note that it is common in the data-driven DRO literature to
suggest choosing $\delta_{n}$ in order to enforce that $P_{\ast}$ is inside
${\mathfrak{M}}_{\delta_{n}}(P_{n})$ with high probability. Such selection
typically \textit{will fall in the third case}, that is, this choice will
induce estimates that are substantially larger than standard statistical
noise. Therefore, prescriptions corresponding to the third case should be
assigned only if the optimizer perceives that the out-of-sample environment is
substantially different from the observed (empirical) environment due to
errors or fluctuations that fall outside of standard statistical noise.

The rest of the paper is organized as follows. In Section \ref{sec-val} we
will quickly review the elements of statistical analysis of Empirical Risk
Minimization (ERM) -- also known as Empirical Optimization or Sample Average
Approximation -- which corresponds to the case $\delta_{n}=0$. Then, in
Section \ref{Sec_DRO_GP}, we will follow a parallel discussion to that of
Section \ref{sec-val} and discuss assumptions for the data-driven DRO version
of the problem. The objective is to use these assumptions so that we can
obtain a flexible and disciplined approach that can be systematically applied
to various DRO formulations. Then, in Section \ref{Section_GP_OV} we will
discuss the application of this approach to the explicit development of
asymptotics for the optimal value in phi-divergence and Wasserstein DRO and,
finally, in Section \ref{Section_GP_OS}, we also develop these explicit
results for associated optimal solutions.

We use the following notation throughout the paper. For a sequence $Y_{n}$ of
random variables, by writing $Y_{n}=o_{p}(n^{-\gamma})$ we mean that
$n^{\gamma}Y_{n}$ tends in probability to zero as $n\to\infty$. In particular
$Y_{n}=o_{p}(1)$ means that $Y_{n}$ tends in probability to zero. The notation
$Q\ll P$ means that $Q\in{\mathfrak{P}}(\mathcal{S})$ is absolutely continuous
with respect to $P\in{\mathfrak{P}}(\mathcal{S})$.
Unless stated otherwise probabilistic statements like ``almost every" (a.e.),  are made with respect to the probability measure $P_{\ast}$.
 By saying that a function
$h:\mathcal{S}\to{\mathbb{R}}$ is integrable we mean that ${\mathbb{E}%
}_{P_{\ast}}|h(X)|<\infty$. It is said that a mapping $\phi:{\mathbb{R}}%
^{m}\to{\mathbb{R}}^{k}$ is directionally differentiable at a point $\theta
\in{\mathbb{R}}^{m}$ if the directional derivative
\begin{equation}
\label{dirdif}\phi^{\prime}(\theta,d):=\lim_{t\downarrow0}\frac{\phi
(\theta+td)-\phi(\theta)}{t}%
\end{equation}
exists for every $d\in{\mathbb{R}}^{m}$. We will use the term $\epsilon
_{n}(\theta)$, $\theta\in\Theta$, to denote a random field such that
\begin{equation}
\label{epsil}\sup_{\theta\in\Theta}\left\vert \epsilon_{n}\left(
\theta\right)  \right\vert =o_{p}(1).
\end{equation}

\setcounter{equation}{0}

\section{Statistics of ERM: Review}

\label{sec-val}

In addition to the population objective function $f(\theta):={\mathbb{E}%
}_{P_{\ast}}\left[  l\left(  X,\theta\right)  \right]  $, introduced in
(\ref{eq-0}), we also let%
\begin{equation}
\vartheta:=\inf_{\theta\in\Theta}f(\theta)\;\;\text{and}\;\;\Theta^{\ast
}:=\mathop{\rm arg\,min}_{\theta\in\Theta}f(\theta), \label{opt-1}%
\end{equation}
be the optimal value and the set of optimal solutions of the population
version of the optimization problem, respectively.

As defined in (\ref{eq-0}), $f_{n}(\theta)={\mathbb{E}}_{P_{n}}\left[
l\left(  X,\theta\right)  \right]  $ is the objective function of the ERM
version of the problem and
\begin{equation}
\vartheta_{n}:=\inf_{\theta\in\Theta}f_{n}(\theta)\;\;\text{and}\;\;\theta
_{n}\in\mathop{\rm arg\,min}_{\theta\in\Theta}f_{n}(\theta) \label{opt-2}%
\end{equation}
are the respective optimal value and an optimal solutions of the ERM problem.
We will now quickly review the development of the asymptotic statistics of the
optimal value in ERM and then we will discuss the corresponding results for
optimal solutions.

\subsection{Asymptotics of the Optimal Value}

\label{sec-optvalue}

In order to analyze the statistical error in the difference between the
optimal values $\vartheta_{n}-\vartheta$, we start from enforcing a functional
Central Limit Theorem (CLT) for $f_{n}\left(  \cdot\right)  $. In particular,
one imposes assumptions which guarantee an expansion of the
form\footnote{Recall that $\epsilon_{n}(\cdot)$ denotes a random field
satisfying condition \eqref{epsil}.}
\begin{equation}
f_{n}\left(  \theta\right)  =f\left(  \theta\right)  +n^{-1/2}r_{n}\left(
\theta\right)  +n^{-1/2}\epsilon_{n}\left(  \theta\right)  ,
\label{eq_FCLT_OF}%
\end{equation}
where we have functional weak convergence
\begin{equation}
r_{n}\left(  \cdot\right)  \dst{\mathfrak{g}} \left(  \cdot\right)
\end{equation}
in the uniform topology on compact sets, with ${\mathfrak{g}}\left(
\cdot\right)  $ being a mean zero Gaussian random field with covariance
function
\begin{equation}
\label{covfun}\mathrm{Cov}\left(  {\mathfrak{g}}\left(  \theta\right)  ,
{\mathfrak{g}}\left(  \theta^{\prime}\right)  \right)  =\mathrm{Cov}_{P_{\ast
}}\big( l\left(  X,\theta\right)  ,l\left(  X,\theta^{\prime}\right)  \big) .
\end{equation}
There are several ways to enforce (\ref{eq_FCLT_OF}); a simple set of
sufficient conditions satisfying this is given next (cf., \cite[example
19.7]{vandervaart}).

\begin{assumption}
\label{ass-1} {\rm (i)} For some $\bar{\theta}\in\Theta$ the expectation
${\mathbb{E}}_{P_{\ast}}[l(X,\bar{\theta})^{2}]$ is finite.  {\rm (ii)}
There is a measurable function $\psi:\mathcal{S}\rightarrow{\mathbb{R}}_{+}$
such that ${\mathbb{E}}_{P_{\ast}}[\psi(X)^{2}]$ is finite and
\begin{equation}
|l(X,\theta)-l(X,\theta^{\prime})|\leq\psi(X)\Vert\theta-\theta^{\prime}%
\Vert\label{lips}%
\end{equation}
for all $\theta,\theta^{\prime}\in\Theta$ and a.e. $X\in\mathcal{S}$.
\end{assumption}

In particular under this assumption, it follows that the expectation function
$f(\theta)$ and variance
\begin{equation}
\label{variance}\sigma^{2}(\theta):=\mathrm{Var}_{P_{\ast}}(l(X,\theta))
\end{equation}
are finite valued and continuous on $\Theta$. Furthermore, since the set
$\Theta$ is compact, it follows that the optimal value $\vartheta_{n}$, of the
ERM problem, converges to $\vartheta$ in probability (in fact almost surely).
Moreover, it is not difficult to show from (\ref{eq_FCLT_OF}) that the
distance from $\theta_{n}$ to $\Theta_{\ast}$ converges in probability to zero
(actually, the convergence occurs almost surely) as $n\rightarrow\infty$.
Finally, since the functional $V(\phi):= \inf_{\theta\in\Theta}\phi(\theta)$,
mapping continuous functions $\phi:\Theta\to{\mathbb{R}}$ to the real line, is
directionally differentiable, the following classical result is a
direct consequence of the (functional) Delta Theorem (cf.,  \cite{shap:91}).

\begin{proposition}
Under Assumption {\rm \ref{ass-1}},
\begin{equation}
n^{1/2}(\vartheta_{n}-\vartheta)\dst\inf_{\theta\in\Theta^{\ast}}%
{\mathfrak{g}}(\theta) \label{dro-a}%
\end{equation}
as $n\rightarrow\infty$. In particular, if $\Theta^{*}=\{\theta^{*}\}$ is the
singleton, i.e. $\theta^{*}$ is the unique optimal solution of the true
problem, then $n^{1/2}(\vartheta_{n}-\vartheta^{*})$ converges in distribution
to normal $N(0,\sigma^{2}(\theta^{*}))$.
\end{proposition}

\subsection{Asymptotics of Optimal Solutions}

\label{sec-optsol}

We assume now that $\Theta^{\ast}=\{\theta^{\ast}\}$ is the singleton, i.e.,
$\theta^{\ast}$ is the unique optimal solution of the true (population)
problem \eqref{eq-2}. We also assume that for a.e. $X$, the function
$l(X,\cdot)$ is continuously differentiable\footnote{Unless stated otherwise
all first and second order derivatives will be taken with respect to vector
$\theta$.}. As it was argued in the previous section, asymptotics of the
optimal value is governed by the asymptotics of the objective function. On the
other hand, asymptotics of the optimal solutions can be derived from the
asymptotics of the gradients of the objective function.

Let us consider the following parametrisation of problem \eqref{eq-2}:
\begin{equation}
\min_{\theta\in\Theta}f\left(  \theta\right)  +v^{T} \theta, \label{eq-2a}%
\end{equation}
with parameter vector $v\in{\mathbb{R}}^{m}$. Denote by $\theta_{*}(v)$ an
optimal solution of the above problem \eqref{eq-2a} viewed as a function of
vector $v$. Of course, we have that $\theta_{*}(0)=\theta^{*}$.

\begin{assumption}
[uniform second order growth]\label{ass-sec} There is a neighborhood
$\mathcal{V}$ of $\theta^{\ast}$ and a positive constant $\kappa$ such that
for every $v$ in a neighborhood of $0\in{\mathbb{R}}^{m}$, problem
\eqref{eq-2a} has an optimal solution $\theta_{*}(v)\in\mathcal{V}$ and
\begin{equation}
f(\theta)+v^{T }(\theta-\theta_{*}\left(  v\right)  )\geq f(\theta_{*}\left(
v\right)  )+\kappa\Vert\theta- \theta_*\left(  v\right)  \Vert^{2},
\label{secorder}%
\end{equation}
for all $\theta\in\Theta\cap\mathcal{V}$.
\end{assumption}

The following assumption can be viewed as a counterpart of Assumption
\ref{ass-1} applied to the gradients of the objective function.

\begin{assumption}
\label{ass-2} {\rm (i)} For some $\bar{\theta}\in\Theta$ the expectation
${\mathbb{E}}_{P_{\ast}}\left[  \Vert\nabla l(X,\bar{\theta})\Vert^{2}\right]
$ is finite.  {\rm (ii)} There is a measurable function $\Psi:\mathcal{S}%
\rightarrow{\mathbb{R}}_{+}$ such that ${\mathbb{E}}_{P_{\ast}}[\Psi(X)^{2}]$
is finite and%
\begin{equation}
\Vert\nabla l(X,\theta)-\nabla l(X,\theta^{\prime})\Vert\leq\Psi(X)\Vert
\theta-\theta^{\prime}\Vert, \label{lips-opt}%
\end{equation}
for all $\theta,\theta^{\prime}\in\Theta$ and a.e. $X\in\mathcal{S}$.
\end{assumption}

By the functional CLT it follows that
\begin{equation}
\nabla f_{n}\left(  \theta\right)  =\nabla f\left(  \theta\right)
+n^{-1/2}d_{n}\left(  \theta\right)  +n^{-1/2}\epsilon_{n}\left(
\theta\right)  , \label{FCLT_GRAD}%
\end{equation}
where we have a functional weak convergence $d_{n}\left(  \cdot\right)
\dst\,{\mathfrak{G}}\left(  \cdot\right)  $ in the uniform topology on a
closed neighborhood of $\theta^{\ast}$, with ${\mathfrak{G}}\left(
\cdot\right)  $ being a continuous mean zero
Gaussian random field with covariance function
\[
\mathrm{Cov}\left[  {\mathfrak{G}}(\theta),{\mathfrak{G}}(\theta^{\prime
})\right]  ={\mathbb{E}}_{P_{\ast}}[(\nabla l(X,\theta)-\nabla f(\theta
))(\nabla l(X,\theta^{\prime})-\nabla f(\theta^{\prime}))^{T}].
\]
It follows from \eqref{FCLT_GRAD} that
\begin{equation}
\left[  \nabla f_{n}(\theta)-\nabla f(\theta)\right]  -\left[  \nabla
f_{n}(\theta^{\ast})-\nabla f(\theta^{\ast})\right]  =n^{-1/2}\left[
d_{n}\left(  \theta\right)  -d_{n}\left(  \theta^{\ast}\right)  +\epsilon
_{n}\left(  \theta\right)  -\epsilon_{n}\left(  \theta^{\ast}\right)  \right]
. \label{delta}%
\end{equation}
Also since $\rho_{n}:=\Vert\theta_{n}-\theta^{\ast}\Vert$ tends in probability
to zero, we have
\begin{equation}
\sup_{\theta:\Vert\theta-\theta^{\ast}\Vert\leq\rho_{n}}\left[  d_{n}\left(
\theta\right)  -d_{n}\left(  \theta^{\ast}\right)  +\epsilon_{n}\left(
\theta\right)  -\epsilon_{n}\left(  \theta^{\ast}\right)  \right]  =o_{p}(1).
\label{delta-2}%
\end{equation}
Thus we have the following result from \cite[Theorem 2.1]{shapmor93}, where
the respective regularity conditions are ensured by the above property \eqref{delta-2}.

\begin{proposition}
\label{pr-equiv} Suppose that Assumptions {\rm \ref{ass-sec}}  and
{\rm \ref{ass-2}}  hold. Then it follows that
\begin{equation}
\label{solas}\theta_{n}=\theta_{*}(Z_{n})+o_{p}(n^{-1/2}),
\end{equation}
where $Z_{n}:=\nabla f_{n}(\theta^{*})-\nabla f(\theta^{*})$.
\end{proposition}

The above result reduces the analysis of asymptotic properties of the optimal
solutions to investigation of asymptotic behavior of the optimal solutions of
the finite dimensional problem \eqref{eq-2a}. By the (finite dimensional)
Central Limit Theorem, $n^{1/2}Z_{n}$ converges in distribution to normal
$N(0,\Sigma)$ with covariance matrix $\Sigma=\mathrm{Cov}(\nabla
l(X,\theta^{*}))$. Moreover, if the mapping $\theta_{*}(v)$ is directionally
differentiable at $v=0$ (in the Hadamard sense), then by the finite
dimensional Delta Theorem it follows from \eqref{solas} that
\begin{equation}
\label{sola-2}n^{1/2}(\theta_{n}-\theta^{*})\,\dst \,\theta_{*}^{\prime}(0,
Z),
\end{equation}
where $Z\sim N(0,\Sigma)$. In particular, if $\theta_{*}^{\prime}(0, w)=A w$
is linear (i.e., $\theta_{*}(v)$ is differentiable at $v=0$ with Jacobian
matrix $A$), then $n^{1/2}(\theta_{n}-\theta^{*})$ converges in distribution
to normal with null mean vector and covariance matrix $A \Sigma A^{T}$.

Directional differentiability of optimal solutions of parameterized problems
is well investigated. For example, if $\theta^{\ast}$ is an interior point of
$\Theta$, $f(\theta)$ is twice continuously differentiable at $\theta^{*}$ and
the Hessian matrix $H:=\nabla^{2}f(\theta^{\ast})$ is nonsingular, then the
uniform second order growth (Assumption \ref{ass-sec}) holds, and
$\theta_{\ast}\left(  v\right)  $ is differentiable at $v=0$ with
$\theta_{\ast}^{\prime}\left(  0,w\right)  =H^{-1}w$. When $\theta^{\ast}$ is
on the boundary of the set $\Theta$, the sensitivity analysis of the
parameterized problem (\ref{eq-2a}) is more delicate and involves a certain
measure of the curvature of the set $\Theta$ at the point $\theta^{\ast}$.
This is discussed extensively in \cite{BS2000}. We also refer to
\cite[sections 5.1.3 and 7.1.5]{SDR} for a basic summary of such results.

It is worthwhile to note at this point that the regularity conditions of
assumptions \ref{ass-sec} and \ref{ass-2} address different properties of the
considered setting. Assumption \ref{ass-sec} deals with the limiting
optimization problem and is of deterministic nature. The \emph{uniform} second
order growth condition was introduced in \cite{shapmor93}, and in a more
general form was discussed in \cite[section 5.1.3]{BS2000}. On the other hand
Assumption \ref{ass-2} is related to the stochastic behavior of the ERM
problem \eqref{eq-1}.

\setcounter{equation}{0}

\section{Statistics of DRO: General Principles}

\label{Sec_DRO_GP}

We now provide sufficient conditions for the development of DRO statistical
principles based on assumptions which are parallel to those imposed in the ERM
section. Define
\begin{equation}
\bar{\vartheta}_{n}:=\inf_{\theta\in\Theta}\mathcal{F}_{n}(\theta,\delta
_{n})\;\;\text{and}\;\;\bar{\theta}_{n}\in\mathop{\rm arg\,min}_{\theta
\in\Theta}\mathcal{F}_{n}(\theta,\delta_{n}), \label{opt-4}%
\end{equation}
the optimal value and an optimal solution of the DRO problem \eqref{eq-3}
(recall the definition of $\mathcal{F}_{n}(\theta,\delta_{n})$ in (\ref{opt-3})).

\subsection{DRO Asymptotics of the Optimal Value}

\label{Subsec_DRO_OV}

Similar to the ERM case, in the DRO setting, we will typically have an
expansion of the form
\begin{equation}
\mathcal{F}_{n}(\theta,\delta_{n})=f_{n}\left(  \theta\right)  +\delta
_{n}^{\gamma}R_{n}\left(  \theta\right)  +\delta_{n}^{\gamma}\epsilon
_{n}\left(  \theta\right)  , \label{DRO_CLT}%
\end{equation}
for some $\gamma>0$, where $R_{n}(\cdot)$ converges in probability in the
uniform topology over $\Theta$ to a continuous deterministic process
$\varrho(\cdot)$,
\begin{equation}
R_{n}(\theta)=\varrho(\theta)+\epsilon_{n}(\theta). \label{DRO_CLT2}%
\end{equation}
Since $\mathcal{F}_{n}(\cdot,\delta_{n})\ge f_{n}(\cdot)$, it follows then
that $\varrho(\cdot)\ge0$.
We will characterize $\gamma>0$ and $\varrho\left(  \cdot\right)  $ explicitly
in the next sections in the context of phi-divergence and Wasserstein DRO
formulations under suitable conditions on the distributional uncertainty set
-- in addition to Assumption \ref{ass-1} (which is clearly independent of the
distributional uncertainty). The following result, which is immediate from the
application of the functional Delta Theorem summarizes the type of behavior
that we expect in DRO\ formulations depending on the geometry of the
distributional uncertainty set and the rate of decay to zero of the
uncertainty size $\delta_{n}$. Recall that ${\mathfrak{g}}(\theta)$ denotes a mean zero Gaussian random field with covariance function \eqref{covfun}.

\begin{theorem}
\label{Thm_GP_DRO_OV} Suppose that Assumption  {\rm \ref{ass-1}} and
conditions   \eqref{DRO_CLT} -  \eqref{DRO_CLT2} hold. Then there are three types of asymptotic
behavior of the DRO optimal value: \newline {\rm (a)} If $\delta_{n}%
^{\gamma}=o\left(  n^{-1/2}\right)  $, then%
\begin{equation}
\bar{\vartheta}_{n}=\vartheta_{n}+o_{p}\left(  n^{-1/2}\right)  ,\label{th-i}%
\end{equation}
and hence
\begin{equation}
n^{1/2}\left(  \bar{\vartheta}_{n}-\vartheta\right)  \dst\inf_{\theta\in
\Theta^{\ast}}{\mathfrak{g}}(\theta),
\end{equation}
which coincides with \eqref{dro-a} and thus the DRO formulation has no
asymptotic impact.
\newline {\rm (b)} If $\delta_{n}^{\gamma}=n^{-1/2}$,
then
\begin{equation}
n^{1/2}\left(  \bar{\vartheta}_{n}-\vartheta\right)  \dst\inf_{\theta\in
\Theta^{\ast}}\left\{  {\mathfrak{g}}(\theta)+\varrho\left(  \theta\right)
\right\}  ,\label{th-ii}%
\end{equation}
so the DRO formulation introduces an explicit and quantifiable asymptotic bias
which can be interpreted as a regularization term.\newline {\rm (c)} If
$o\left(  \delta_{n}^{\gamma}\right)  =n^{-1/2}$, then\footnote{The right hand
side of \eqref{th-iii} is a deterministic number. Therefore convergence in
distribution `$\dst$' there is the same as convergence in probability.}
\begin{equation}
\delta_{n}^{-\gamma}\left(  \bar{\vartheta}_{n}-\vartheta\right)
\dst\inf_{\theta\in\Theta^{\ast}}\varrho\left(  \theta\right)  ,\label{th-iii}%
\end{equation}
so the bias term induced by the DRO\ formulation is larger than the
statistical error.
\end{theorem}

\textbf{Proof.} Part  (a).  By \eqref{DRO_CLT} and
\eqref{DRO_CLT2} we have that in the considered case
$$ \F_n(\theta,\delta_n)=f_n(\theta)+ o(n^{-1/2})\epsilon
_{n}\left(  \theta\right),$$
where   $\epsilon_{n}\left(  \theta\right)$ is the generic term satisfying \eqref{epsil}.
Thus \eqref{th-i} follows.

 Part (b).  By \eqref{DRO_CLT} and
\eqref{DRO_CLT2}  in the considered case we can write
\[
 n^{1/2}\left(\mathcal{F}_{n}(\theta,\delta_{n})-f(\theta)\right)=  n^{1/2}(f_{n} (  \theta  )-f(\theta))  + \varrho (\theta)+ \epsilon
_{n}\left(  \theta\right).
\]
Under Assumption  {\rm \ref{ass-1}},   by the functional CLT
we have that $n^{1/2}(f_{n} (  \theta  )-f(\theta))  + \varrho (\theta)$ converges in distribution to $\cg(\theta)+\varrho(\theta)$,
and hence \eqref{th-ii} follows by the Delta Theorem.

 Part
(c)\ may appear somewhat different because the right hand side is deterministic
but, under case (c)\ note that we can simply write
\[
\mathcal{F}_{n}(\theta,\delta_{n})=f\left(  \theta\right)  +\delta_{n}%
^{\gamma}R_{n}\left(  \theta\right)  +\delta_{n}^{\gamma}\epsilon_{n}\left(
\theta\right)  ,
\]
so case (c) also follows from the standard analysis since $R_{n}\left(
\cdot\right)  $ converges uniformly to $\varrho\left(  \cdot\right)  $ in
probability (thus it converges weakly in the uniform topology). $\hfill\square$

\subsection{DRO Asymptotics of the Optimal Solutions}

\label{Subsec_DRO_OS}

As in the ERM development, in addition to Assumptions \ref{ass-sec},
it is convenient to guarantee that for all $n$ large enough,
$\mathcal{F}_{n}(\theta,\delta_{n})$ is differentiable in a neighborhood
$\mathcal{V}$ of $\theta^{\ast}$ and%
\begin{equation}
\nabla\mathcal{F}_{n}(\theta,\delta_{n})=\nabla f_{n}\left(  \theta\right)
+\delta_{n}^{\gamma}D_{n}\left(  \theta\right)  +\delta_{n}^{\gamma}%
\epsilon_{n}\left(  \theta\right)  , \label{DROGCLT}%
\end{equation}
for some $\gamma>0$, where $D_{n}\left(  \theta\right)  $ converges in
probability to $\nabla\varrho\left(  \theta\right)  $ uniformly around a
closed neighborhood $\mathcal{\bar{V}}$ of $\theta^{\ast}$.
In consequence, we obtain the following analog of Theorem \ref{Thm_GP_DRO_OV},
which follows from the \textit{finite dimensional} Delta Theorem. Recall that
$\theta_{*}(v)$ is an optimal solution of problem \eqref{eq-2a} and
$\theta^{\prime}_{*}(0,\cdot)$ is its directional derivative at $v=0$.

\begin{theorem}
\label{Thm_GP_DRO_OS} Suppose that: Assumptions {\rm \ref{ass-sec}
} and {\rm \ref{ass-2}} hold, conditions   \eqref{DRO_CLT} -  \eqref{DRO_CLT2} are satisfied,
 identity \eqref{DROGCLT} holds with
$D_{n}\left(  \cdot\right)  $ converging in probability to $\nabla
\varrho\left(  \cdot\right)  $ uniformly around a closed neighborhood
$\mathcal{\bar{V}}$ of $\theta^{\ast}$, and that $\theta_{\ast}(v)$ is
directionally differentiable at $v=0$ (in the Hadamard sense). Let $Z\sim
N(0,\Sigma)$ with covariance matrix $\Sigma=\mathrm{Cov}(\nabla l(X,\theta
^{\ast}))$. Then the DRO optimal solutions can have three types of asymptotic
behavior: \newline {\rm (A)} If $\delta_{n}^{\gamma}=o\left(  n^{-1/2}%
\right)  $, then%
\begin{equation}
\bar{\theta}_{n}=\theta_{n}+o_{p}\left(  n^{-1/2}\right)  ,
\end{equation}
thus
\begin{equation}
n^{1/2}\left(  \bar{\theta}_{n}-\theta^{\ast}\right)  \dst\,\theta_{\ast
}^{\prime}\left(  0,Z\right)  \mathfrak{.}%
\end{equation}
\newline {\rm (B)} If $\delta_{n}^{\gamma}=n^{-1/2}$, then
\begin{equation}
n^{1/2}\left(  \bar{\theta}_{n}-\theta^{\ast}\right)  \dst\,\theta_{\ast
}^{\prime}\left(  0,Z+\nabla\varrho(\theta^{\ast})\right)  .
\end{equation}
\newline {\rm (C)} If $o\left(  \delta_{n}^{\gamma}\right)  =n^{-1/2}$,
then
\begin{equation}
\delta_{n}^{-\gamma}\left(  \bar{\theta}_{n}-\theta^{\ast}\right)
\dst\,\theta_{\ast}^{\prime}\left(  0,\nabla\varrho\left(  \theta^{\ast
}\right)  \right)  .
\end{equation}

\end{theorem}

\setcounter{equation}{0}

\section{General Principle in Action: Optimal Values}

\label{Section_GP_OV}

In this section, we apply the general principle to the asymptotics of the
value function in two of the main types of DRO formulations, namely,
phi-divergence and Wasserstein DRO.

\subsection{The Phi-Divergence Case}

We recall the definition of the distributional uncertainty set for the
phi-divergence case. Consider a convex lower semi-continuous function
$\phi\colon{\mathbb{R}}\rightarrow{\mathbb{R}}_{+}\cup\{+\infty\}$ such that
$\phi(1)=0$ and $\phi(t)=+\infty$ for $t<0$. For probability measures
$Q,P\in{\mathfrak{P}}$ such that $Q$ is absolutely continuous with respect to
$P$ with the corresponding density $dQ/dP$, the $\phi$-divergence is defined
as (cf., \cite{Csis1963},\cite{Morimoto1963})
\begin{equation}
D_{\phi}(Q\Vert P):={\mathbb{E}}_{P}[\phi(dQ/dP)]=\int\phi(dQ/dP)dP.
\label{phidiver}%
\end{equation}
In particular, for $\phi(t):=t\log\left(  t\right)  -t+1$, $t\geq0$, this
becomes the Kullback--Leibler (KL) divergence of~$Q$ from~$P$. The ambiguity
set ${\mathfrak{M}}_{\delta}(P)$ associated with $D_{\phi}(\cdot\Vert P)$ is
defined as
\begin{equation}
{\mathfrak{M}}_{\delta}\left(  P\right)  :=\{Q\ll P:D_{\phi}(Q\Vert
P)\leq\delta\}. \label{phiset}%
\end{equation}
By duality arguments the corresponding distributionally robust functional can
be written in the form (cf., \cite{bl-2015}, \cite{bental}, \cite{sha2017})
\begin{equation}
\sup_{Q\in{\mathfrak{M}}_{\delta}\left(  P\right)  }{\mathbb{E}}_{Q}%
[Y]=\inf_{\mu,\lambda>0}\left\{  \lambda\delta+\mu+\lambda{\mathbb{E}}%
_{P}[\phi^{\ast}((Y-\mu)/\lambda)]\right\}  , \label{diverg-1}%
\end{equation}
where $\phi^{\ast}(y)=\sup_{t\in{\mathbb{R}}}\{yt-\phi(t)\}$ is the convex
conjugate of $\phi$. Using this representation we can obtain an asymptotic
expansion for (\ref{diverg-1}) as a function of $\delta$. This expansion can
be helpful to suggest the form of the expansion in (\ref{DRO_CLT}) and
(\ref{DROGCLT}). For this, we need to assume certain regularity properties of
$\phi\left(  t\right)  $ around $t=1$.

\begin{assumption}
\label{ass-phi} Assume that $\phi(t)$ is two times continuously differentiable
in a neighborhood of $t=1$ with $\kappa:=2/\phi^{\prime\prime}(1)>0$.
\end{assumption}

Under this condition we have the following expansion, which is obtained, in
order to simplify our exposition, under the assumption that the probability
measure $P$ has compact support. See also the results in \cite{Lam16}, which provide additional correction terms under a fixed $P$. The uniform feature of the statement below is helpful in the statistical analysis. Our development here will also be used in the expansion of the optimal solutions.

\begin{proposition}
\label{Prop_Expand_Phi}Suppose that Assumption  {\rm \ref{ass-phi}} holds,
that $P\left(  \left\vert Y\right\vert \leq\nu\right)  =1$ for some $\nu
\in\left(  0,\infty\right)  $. Then, for any $b_{0}>0$,%
\begin{equation}
\sup_{Q\in{\mathfrak{M}}_{\delta}\left(  P\right)  }{\mathbb{E}}%
_{Q}[Y]-{\mathbb{E}}_{P}(Y)-\delta^{1/2}\kappa^{1/2}\sqrt{\mathrm{Var}_{P}%
[Y]}=o\left(  \delta^{1/2}\right)  , \label{As_delta}%
\end{equation}
uniformly over
Borel probability measures $P$ supported on $[-\nu,\nu]$  such that
$\mathrm{Var}_{P}[Y]\geq b_{0}$. Moreover, there is $\bar{\delta}>0$ such that
for all $\delta<\bar{\delta}$
\[
\arg\max\{{\mathbb{E}}_{Q}[Y]:Q\in{\mathfrak{M}}_{\delta}\left(  P\right)  \}
\]
is unique.
\end{proposition}

\textbf{Proof.} Note that we can write
\[
\sup_{Q\in{\mathfrak{M}}_{\delta}\left(  P\right)  }{\mathbb{E}}_{Q}%
[Y]=\sup_{{\mathbb{E}}_{P}\left(  Z\right)  =1,{\mathbb{E}}_{P}\left(
\phi\left(  Z\right)  \right)  \leq\delta}{\mathbb{E}}_{P}[YZ],
\]
where the sup is taken over the set of positive random variables $Z$
satisfying the specified moment constraints. We may assume that ${\mathbb{E}%
}_{P}\left[  Y\right]  =0$ for simplicity since we can always center the
objective function around ${\mathbb{E}}_{P}\left[  Y\right)  ]$. In turn, by
letting $\bar{\Delta}=(Z-1)/\delta^{1/2}$, the previous optimization problem
is equivalent to
\begin{equation}
\delta^{1/2}\sup_{{\mathbb{E}}_{P}\left(  \bar{\Delta}\right)  =0,\bar{\Delta
}\geq-\delta^{-1/2},{\mathbb{E}}_{P}(\phi(1+\delta^{1/2}\bar{\Delta}%
))\leq\delta}{\mathbb{E}}_{P}[Y\bar{\Delta}]. \label{OPT_PHI}%
\end{equation}
Since $\left\vert Y\right\vert \leq\nu$ and ${\mathbb{E}}_{P}\left[  Y\right]
=0$, then $\bar{\Delta}=aY$ is feasible for any $a>0$ provided that $a \nu
\leq\delta^{-1/2}$ and
\[
{\mathbb{E}}_{P}[\phi(1+\delta^{1/2}\bar{\Delta})]\leq\delta.
\]
In turn, since $\phi\left(  t\right)  $ is two times continuously
differentiable at $t=1$, we have that
\[
\delta^{-1}\phi(1+\delta^{1/2}ay))\rightarrow a^{2}y^{2}\phi^{\prime\prime
}\left(  1\right)  /2
\]
as $\delta\rightarrow0$ uniformly over compact sets. Therefore, we conclude
that there exists $\delta_{0}>0$ such that for any $\delta<\delta_{0}$
\begin{align*}
&  \sup_{{\mathbb{E}}_{P}\left(  \bar{\Delta}\right)  =0,\bar{\Delta}%
\geq-\delta^{-1/2},{\mathbb{E}}_{P}(\phi(1+\delta^{1/2}\bar{\Delta}%
))\leq\delta}{\mathbb{E}}_{P}[Y\bar{\Delta}]\\
&  \geq\sup_{a>0,a^{2}{\mathbb{E}}_{P}\left(  Y^{2}\right)  /2\leq\left(
1-\delta_{0}\right)  /\phi^{\prime\prime}\left(  1\right)  }{\mathbb{E}}%
_{P}[aY^{2}]=\sqrt{\kappa\left(  1-\delta_{0}\right)  }\cdot\sqrt{{\mathbb{E}%
}_{P}[Y^{2}]}.
\end{align*}
Since $\delta_{0}>0$ can be chosen to be arbitrarily small, we conclude an
asymptotic lower bound which retrieves (\ref{As_delta}). For the upper bound,
we apply the duality result (\ref{diverg-1}) in the form corresponding to
(\ref{OPT_PHI}), we obtain
\begin{align}
&  \sup_{{\mathbb{E}}_{P}\left(  \bar{\Delta}\right)  =0,\bar{\Delta}%
\geq-\delta^{-1/2},\delta^{-1}{\mathbb{E}}_{P}(\phi(1+\delta^{1/2}\bar{\Delta
}))\leq1}{\mathbb{E}}_{P}[Y\bar{\Delta}]\nonumber\\
&  =\min_{\bar{\lambda}>0,\bar{\mu}}\{\bar{\lambda}+{\mathbb{E}}_{P}%
[\sup_{\bar{\Delta}\geq-\delta^{-1/2}}\{\left(  Y+\bar{\mu}\right)
\bar{\Delta}-\bar{\lambda}\delta^{-1/2}\phi(1+\delta^{1/2}\bar{\Delta
})\}]\}\nonumber\\
&  \leq\min_{\bar{\lambda}>0}\{\bar{\lambda}+{\mathbb{E}}_{P}[\sup
_{\bar{\Delta}\geq-\delta^{-1/2}}\{Y\bar{\Delta}-\bar{\lambda}\delta
^{-1/2}\phi(1+\delta^{1/2}\bar{\Delta})\}]\}. \label{Stab_PhiDiv}%
\end{align}
We will plug in
\[
\bar{\lambda}_{0}=\arg\min\{\bar{\lambda}+\kappa{\mathbb{E}}_{P}[Y^{2}%
]/4\bar{\lambda}:\bar{\lambda}>0\}=2^{-1}\sqrt{\kappa{\mathbb{E}}_{P}[Y^{2}%
]}>0
\]
into (\ref{Stab_PhiDiv}) to obtain our upper bound. Using that $\bar{\lambda
}_{0}>0$ and that $\phi$ is convex with $\phi^{\prime\prime}\left(  1\right)
>0$, we have that the family of (continuous) functions
\[
s_{\delta}\left(  y\right)  :=\sup_{\bar{\Delta}\geq-\delta^{-1/2}}%
\{y\bar{\Delta}-\bar{\lambda}\delta^{-1/2}\phi(1+\delta^{1/2}\bar{\Delta})\}
\]
converges uniformly on compact sets to
\[
s_{0}\left(  y\right)  =\sup_{\bar{\Delta}}\{y\bar{\Delta}-\bar{\lambda}%
\bar{\Delta}^{2}/\kappa\}=\frac{\kappa y^{2}}{4\bar{\lambda}}.
\]
Therefore we obtain that
\begin{align*}
&  \min_{\bar{\lambda}>0}\{\bar{\lambda}+{\mathbb{E}}_{P}[\sup_{\bar{\Delta
}\geq-\delta^{-1/2}}\{Y\bar{\Delta}-\bar{\lambda}\delta^{-1/2}\phi
(1+\delta^{1/2}\bar{\Delta})\}]\}\\
&  \leq\bar{\lambda}_{0}+{\mathbb{E}}_{P}[\sup_{\bar{\Delta}\geq-\delta
^{-1/2}}\{Y\bar{\Delta}-\bar{\lambda}_{0}\delta^{-1/2}\phi(1+\delta^{1/2}%
\bar{\Delta})\}]\rightarrow\sqrt{\kappa}\cdot\sqrt{{\mathbb{E}}_{P}[Y^{2}]}.
\end{align*}
These estimates, which are uniform given that $\left\vert Y\right\vert \leq
\nu$, yield the estimate in the proposition. The uniqueness is standard, it
follows from the local strong convexity of $\phi\left(  \cdot\right)  $ at the
origin.
$\hfill\square$ \newline

Recall that  $\sigma^{2}(\theta):=\mathrm{Var}_{P_{\ast}}(l(X,\theta
))$, and that
${\mathfrak{g}}\left(
\cdot\right) $  is  a mean zero Gaussian random field.
 Expansion (\ref{As_delta}) immediately yields, at least when $\sup_{\theta
\in\Theta}\left\vert l\left(  X,\theta\right)  \right\vert $ is $P_{\ast}%
$-bounded, that
\begin{equation}
\mathcal{F}_{n}(\theta,\delta_{n})=f_{n}\left(  \theta\right)  +\delta
_{n}^{1/2}\kappa^{1/2} \sigma(\theta) +\delta_{n}^{1/2}\epsilon_{n}\left(
\theta\right)  . \label{empest-1}%
\end{equation}
  Consequently, we obtain the following result.

\begin{theorem}
\label{th-asymdro} Suppose that $\sup_{\theta\in\Theta}\left\vert l\left(
X,\theta\right)  \right\vert $ is $P_{\ast}$-essentially bounded, that
Assumption  {\rm \ref{ass-1}} and Assumption  {\rm \ref{ass-phi}} hold,
and that $\sigma^{2}(\theta)>0$ for all $\theta\in\Theta^{*}$. Then, we have
the following types of asymptotic behavior of the DRO optimal values.\newline%
 {\rm (a-phi)} If $\delta_{n}=o\left(  n^{-1}\right)  $, then%
\begin{equation}
n^{1/2}\left(  \bar{\vartheta}_{n}-\vartheta\right)  \dst\inf_{\theta\in
\Theta^{\ast}}{\mathfrak{g}}(\theta).
\end{equation}
 {\rm (b-phi)} If $\delta_{n}=\beta n^{-1}$, then
\begin{equation}
\label{b-phi}n^{1/2}\left(  \bar{\vartheta}_{n}-\vartheta\right)
\dst\inf_{\theta\in\Theta^{\ast}}\left\{  {\mathfrak{g}}(\theta)+ \kappa^{1/2}
\beta^{1/2} \sigma(\theta)\right\}  .
\end{equation}
 {\rm (c-phi)} If $o\left(  \delta_{n}\right)  =n^{-1}$, then
\begin{equation}
\label{c-phi}\delta_{n}^{-1/2}\left(  \bar{\vartheta}_{n}-\vartheta\right)
\dst\,\kappa^{1/2}\inf_{\theta\in\Theta^{\ast}}\sigma(\theta),
\end{equation}
so the bias term induced by the DRO\ formulation dominates the statistical error.
\end{theorem}

\textbf{Proof.} Proof of this theorem is quite standard (cf., \cite[ proof of
Theorem 5.7]{SDR}). For the sake of completeness we briefly outline proof of
case (b-phi). Note that our assumptions imply Assumption \ref{ass-1}, and
hence $\sigma^{2}(\theta)$ is a continuous function of $\theta$. Therefore
there is a compact neighborhood $\bar{\Theta}$ of $\Theta^{*}$ such that
$\sigma^{2}(\theta)>0$ for all $\theta\in\bar{\Theta}$. We can restrict the
minimization to $\bar{\Theta}$ for which the expansion \eqref{empest-1} holds.

Consider the space $C(\bar{\Theta})$ of continuous functions $g:\bar{\Theta
}\to{\mathbb{R}}$ equipped with the sup-norm, and functional $V(g):=\inf
_{\theta\in\bar{\Theta}} g(\theta)$, mapping $C(\bar{\Theta})$ into the real
line. This functional is directionally differentiable in the Hadamard sense
with the directional derivative at a point $\mu\in C(\bar{\Theta})$ given by
$V^{\prime}(\mu, h)=\inf_{\theta\in\bar{\Theta}(\mu)} h(\theta)$, where
$\bar{\Theta}(\mu):=\mathop{\rm arg\,min}_{\theta\in\bar{\Theta}} \mu(\theta
)$. We have that $\bar{\vartheta}_{n}=V(\mathcal{F}_{n})$ and $\vartheta
=V(f)$, where $\mathcal{F}_{n}(\cdot):=\mathcal{F}_{n}(\cdot,\delta_{n})$. By
the functional CLT and \eqref{empest-1} it follows that $n^{1/2}%
(\mathcal{F}_{n}-f)$ converges in distribution (weakly) to ${\mathfrak{g}%
}(\theta)+ \kappa^{1/2}\beta^{1/2}\sigma(\theta)$. We can apply now the
functional Delta Theorem to conclude \eqref{b-phi}. $\hfill\square$ \newline


Given that $\phi(\cdot)$ is only assumed to satisfy Assumption \ref{ass-phi},
without imposing any growth condition, situations such as the (c-phi) case
require imposing stronger moment conditions than just assuming $\mathrm{Var}%
_{P_{\ast}}[l\left(  X,\theta\right)  ]<\infty$. This can be seen in the
KL-divergence case in which $\phi\left(  t\right)  =t\log\left(  t\right)
-t+1$. For fixed $\delta>0$, the population solution requires that $l\left(
X,\theta\right)  $ has a finite moment generating function in a neighborhood
of the origin. Therefore, if $\delta_{n}$ converges to zero sufficiently
slowly and $l\left(  X,\theta\right)  $ has infinite moments of order
$2+\varepsilon$, an expansion such as (\ref{empest-1}) may not hold. However,
if $\phi\left(  t\right)  =\left(  t-1\right)  ^{2}$, it follows that
expansion (\ref{empest-1}) holds exactly with $\epsilon_{n}\left(
\theta\right)  =0$.

On the other hand, the result in \cite[Theorem 2]{Duchi:2021} provides more
for the case (b-phi) since it does not require compact support (although it
requires $\phi$ to be three times continuously differentiable).\ The following
example shows that the smoothness of $\phi\left(  \cdot\right)  $ is important
in deriving the asymptotics in the previous result with $\delta_{n}=n^{-1/2}$.

\begin{example}
\label{ex-abs}  {\rm Consider $\phi(t):=|t-1|$, $t\geq0$. In that case
(e.g., \cite[Example 3.12]{sha2017}), for $\delta\in(0,2)$ and essentially
bounded $Y$,
\begin{equation}
\sup_{Q\in{\mathfrak{M}}_{\delta}(P)}{\mathbb{E}}_{Q}[Y]=(\delta
/2)\ess(Y)+(1-\delta/2)\mathsf{AV@R}_{P,1-\delta/2}(Y),\label{abs-1}%
\end{equation}
where
\begin{equation}
\mathsf{AV@R}_{P,\alpha}(Y):=\inf_{\tau\in{\mathbb{R}}}\left\{  \tau
+\alpha^{-1}{\mathbb{E}}_{P}[Y-\tau]_{+}\right\}  ,\;\alpha\in
(0,1].\label{abs-2}%
\end{equation}
Note that $\mathsf{AV@R}_{P,1}(Y)={\mathbb{E}}_{P}[Y]$ and as $\alpha$ tends
to one,
\begin{equation}
\big|\mathsf{AV@R}_{P,\alpha}(Y)-{\mathbb{E}}_{P}[Y]\big|=O(1-\alpha
),\label{abs-3}%
\end{equation}
provided $Y$ is essentially bounded.

 Suppose that $l(x,\theta)$ is bounded on $\mathcal{S}\times\Theta$,
and hence
\begin{equation}
\mathcal{F}_{n}(\theta,\delta_{n})=(\delta_{n}/2)\max_{1\leq i\leq N}%
l(X_{i},\theta)+(1-\delta_{n}/2)\mathsf{AV@R}_{P_{n},1-\delta_{n}%
/2}(l(X,\theta)).\label{abs-4}%
\end{equation}
Consider $\mathrm{\delta}$ $_{n}=\beta n^{-1}$ with $\beta>0$. Then the first
term in \eqref{abs-4} is of order $O(n^{-1})$, and by \eqref{abs-3} the second
term is ${\mathbb{E}}_{P_{n}}[l(X,\theta)]+O\left(  n^{-1}\right)  $.
Consequently in that case $\bar{\vartheta}_{n}=\vartheta_{n}+o_{p}(n^{-1/2}),$
and hence this corresponds to case (a)\ in Theorem \ref{Thm_GP_DRO_OV}. This
shows that the assumption of smoothness (differentiability) of $\phi(\cdot)$
is essential for derivation of the asymptotics of $\bar{\vartheta}_{n}$. Here
some additional terms in the asymptotics of $\bar{\vartheta}_{n}$ appear when
$\delta_{n}$ is of order $O(n^{-1/2})$, rather than $O(n^{-1})$.
$\hfill\square$ }
\end{example}


\subsection{The Wasserstein Distance Case}

We use ${\mathfrak{P}}(\mathcal{S}\times\mathcal{S})$ to denote the set of
Borel probability measures on the product space $\mathcal{S\times S}$. Let
$c:\mathcal{S\times S}\rightarrow{\mathbb{R}}_{+}\cup\{+\infty\}$ be a lower
semi-continuous function such that $c(x,y)=0$ if and only if $x=y$. This
function measures the marginal cost of a transporting a unit of mass from a
source location to a target location, respectively.
The optimal transport cost between $P,Q\in{\mathfrak{P}}(\mathcal{S})$ is
given by
\begin{equation}
D_{c}\left(  P,Q\right)  :=\min\{{\mathbb{E}}_{\pi}\left[  c\left(
X,Y\right)  \right]  :\pi\in{\mathfrak{P}}\left(  \mathcal{S}\times
\mathcal{S}\right)  ,\text{ }\pi_{X}=P,\text{ }\pi_{Y}=Q\}, \label{cost}%
\end{equation}
where ${\mathbb{E}}_{\pi}[\,\cdot\,]$ is the expectation under a joint
distribution $\pi\in{\mathfrak{P}}\left(  \mathcal{S}\times\mathcal{S}\right)
$ and $\pi_{X}$ and $\pi_{Y}$ denote the marginal distributions of $X$ and
$Y$, respectively. It turns out that the optimizer is always achieved, thus we
write `$\min$' instead of $`\inf$'. Let $\Vert\cdot\Vert$ be a norm on the
space ${\mathbb{R}}^{d}$. An important special case corresponds to the choice
$c\left(  x,y\right)  :=\Vert x-y\Vert^{p}$ for some $p>0$, in which case
$D_{c}\left(  P,Q\right)  ^{1/p}$ is the so-called $p$-Wasserstein distance.
The reader is referred to the text of Villani \cite{villani} for more
background on optimal transport.

For any given $P\in{\mathfrak{P}}(\mathcal{S})$ and $\delta\geq0$ we have the
following dual result (cf., \cite{esf-2015}, \cite{BKM-2019}, \cite{GK16}) assuming that ${\mathfrak{h}%
}\left(  \cdot\right)  $ is upper semi-continuous and ${\mathfrak{h}}(X)$ is
$P$-integrable,%
\begin{equation}
\sup_{Q:\,D_{c}(P,Q)\leq\delta}{\mathbb{E}}_{Q}[{\mathfrak{h}}(Y)]=\min
_{\lambda\geq0}\left\{  \lambda\delta+{\mathbb{E}}_{P}\left[  {\mathfrak{\bar
{h}}}_{\lambda}(X)\right]  \right\}  , \label{wasdual-1}%
\end{equation}
where
\begin{equation}
{\mathfrak{\bar{h}}}_{\lambda}(x):=\sup_{y\in\mathcal{S}}\{{\mathfrak{h}%
}(y)-\lambda c(x,y)\},\;\lambda\ge0. \label{wasdual-2}%
\end{equation}

Throughout the rest of our discussion, we will choose $c\left(  x,y\right)
:=\Vert x-y\Vert^{p}$ for $p\in\left(  1,\infty\right)  $ and therefore write
$D_{p}\left(  P,Q\right)  $ for this choice of cost function. Further, we use
$\left\Vert \cdot\right\Vert _{\ast}$ to denote the dual norm, namely,
\[
\left\Vert y\right\Vert _{\ast}=\sup\{x^{T}y:\left\Vert x\right\Vert =1\}.
\]

As in the phi-divergence case, assuming that $P$ is fixed and has compact
support, for example, we can obtain an asymptotic expansion for
(\ref{wasdual-1}) as a function of $\delta$.  By writing $\bbe_P^{(p-1)/p}[\,\cdot\,]$ we mean $(\bbe_P[\,\cdot\,])^{(p-1)/p}$.

\begin{proposition}
\label{Prof_Sens_WDRO}Suppose that ${\mathfrak{h}}\left(  \cdot\right)  $ is
continuously differentiable and the mapping
\begin{equation}
x\mapsto\sup\{\left\Vert \nabla{\mathfrak{h}}\left(  x+\Delta\right)
-\nabla{\mathfrak{h}}\left(  x\right)  \right\Vert /(1+\left\Vert
\Delta\right\Vert ^{p-1}):\Delta\in{\mathbb{R}}^{d}\} \label{Loc_bdd}%
\end{equation}
is bounded on compact sets. Then, for any  $b_{0}>0$,
\[
\sup_{Q:\,D_{p}(P,Q)\leq\delta}{\mathbb{E}}_{Q}[{\mathfrak{h}}(Y)]-{\mathbb{E}%
}_{P}[{\mathfrak{h}}(X)]-\delta^{1/p}{\mathbb{E}}_{P}^{\left(  p-1\right)
/p}[\left\Vert \nabla{\mathfrak{h}}(X)\right\Vert _{\ast}^{p/\left(
p-1\right)  }]=o\left(  \delta^{1/p}\right)  ,
\]
uniformly over $P\in{\mathfrak{P}}(\mathcal{[}-\nu,\nu\mathcal{]}^{d})$ such that 
${\mathbb{E}}_{P}\left\Vert \nabla{\mathfrak{h}}\left(  Y\right)  \right\Vert
\geq b_{0}$.
\end{proposition}

\textbf{Proof.} The proof of this result is similar to the one given in the
phi-divergence case. We start by observing that
\[
\sup_{Q:\,D_{p}(P,Q)\leq\delta}{\mathbb{E}}_{Q}[{\mathfrak{h}}(Y)]={\mathbb{E}%
}_{P}[{\mathfrak{h}}(X)]+\sup_{{\mathbb{E}}_{P}\left\Vert \Delta\right\Vert
^{p}\leq\delta}{\mathbb{E}}_{P}[{\mathfrak{h}}(X+\Delta)-{\mathfrak{h}}(X)],
\]
where the optimization in the right hand side is taken over random variables
$\Delta$. We let $\delta^{1/p}\bar{\Delta}=\Delta$ and note that
\begin{align*}
&  \sup_{{\mathbb{E}}_{P}\left\Vert \Delta\right\Vert ^{p}\leq\delta
}{\mathbb{E}}_{P}[{\mathfrak{h}}(X+\Delta)-{\mathfrak{h}}(X)]\\
&  =\delta^{1/p}\sup_{{\mathbb{E}}_{P}\left\Vert \bar{\Delta}\right\Vert
^{p}\leq1}{\mathbb{E}}_{P}[\left(  {\mathfrak{h}}(X+\delta^{1/p}\bar{\Delta
})-{\mathfrak{h}}(X)\right)  /\delta^{1/p}]\\
&  =\delta^{1/p}\sup_{{\mathbb{E}}_{P}\left\Vert \bar{\Delta}\right\Vert
^{p}\leq1}{\mathbb{E}}_{P}\left[  \int_{0}^{1}\nabla{\mathfrak{h}}%
(X+t\delta^{1/p}\bar{\Delta})\cdot\bar{\Delta}dt\right]  .
\end{align*}
Next, we can obtain a lower bound by considering a specific form of
$\bar{\Delta}$ suggested by the formal asymptotic limit as $\delta
\rightarrow0$. Note that
\[
{\mathbb{E}}_{P}[\nabla{\mathfrak{h}}(X)\cdot\bar{\Delta}]\leq{\mathbb{E}}%
_{P}[\left\Vert \nabla{\mathfrak{h}}(X)\right\Vert _{\ast}\left\Vert
\bar{\Delta}\right\Vert ],
\]
and the equality is achieved if we choose any $\bar{\Delta}_{\ast}$ which is a
constant multiple of
\[
\bar{\Delta}_{1}\left(  X\right)  \in\arg\max\{\nabla{\mathfrak{h}}%
(X)\cdot\bar{\Delta}:\left\Vert \bar{\Delta}\right\Vert =1\},
\]
(The function $\bar{\Delta}_{1}\left(  \cdot\right)  $ can be selected in a
measurable way using the uniformization technique of Jankov-von Neumann.)
Next, if $\left\Vert \bar{\Delta}^{\ast}\right\Vert =a\left\Vert
\nabla{\mathfrak{h}}(X)\right\Vert _{\ast}^{\gamma}$, then
\[
{\mathbb{E}}_{P}[\left\Vert \nabla{\mathfrak{h}}(X)\right\Vert _{\ast
}\left\Vert \bar{\Delta}^{\ast}\right\Vert ]=a{\mathbb{E}}_{P}[\left\Vert
\nabla{\mathfrak{h}}(X)\right\Vert _{\ast}^{\gamma+1}]
\]
and
\[
{\mathbb{E}}_{P}\left(  \left\Vert \bar{\Delta}^{\ast}\right\Vert ^{p}\right)
=a^{p}{\mathbb{E}}_{P}\left\Vert \nabla{\mathfrak{h}}(X)\right\Vert _{\ast
}^{\gamma p}=1.
\]
Letting $\gamma p=\gamma+1$ we have that $\gamma=1/(p-1)$ and therefore
\[
\sup_{{\mathbb{E}}_{P}\left\Vert \bar{\Delta}\right\Vert ^{p}\leq1}%
{\mathbb{E}}_{P}[\nabla{\mathfrak{h}}(X)\cdot\bar{\Delta}^{\ast}]={\mathbb{E}%
}_{P}^{\left(  p-1\right)  /p}[\left\Vert \nabla{\mathfrak{h}}(X)\right\Vert
_{\ast}^{p/\left(  p-1\right)  }],
\]
with
\[
\bar{\Delta}^{\ast}\left(  X\right)  =\bar{\Delta}_{1}\left(  X\right)
\left\Vert \nabla{\mathfrak{h}}(X)\right\Vert _{\ast}^{1/\left(  p-1\right)
}{\mathbb{E}}_{P}^{-1/p}\left\Vert \nabla{\mathfrak{h}}(X)\right\Vert _{\ast
}^{p/\left(  p-1\right)  }.
\]
The denominator is well defined since ${\mathbb{E}}_{P}\left\Vert
\nabla{\mathfrak{h}}\left(  Y\right)  \right\Vert >0$ and the random variable
$\bar{\Delta}^{\ast}\left(  X\right)  $ is essentially bounded uniformly over
the family $P\in{\mathfrak{P}}(\mathcal{[}-\nu,\nu\mathcal{]}^{d})$ and
${\mathbb{E}}_{P}\left\Vert \nabla{\mathfrak{h}}\left(  Y\right)  \right\Vert
\geq b_{0}$. Since the gradient $\nabla{\mathfrak{h}}(\cdot)$ is continuous,
then it is uniformly continuous over compact sets and, consequently, uniformly
over $\bar{\Delta}$ in compact sets,%
\[
\int_{0}^{1}\left\Vert \nabla{\mathfrak{h}}(x+t\delta^{1/p}\bar{\Delta
})-\nabla{\mathfrak{h}}(x)\right\Vert \bar{\Delta}dt=o\left(  1\right)
\]
as $\delta\rightarrow0$. This yields that
\[
\sup_{{\mathbb{E}}_{P}\left\Vert \bar{\Delta}\right\Vert ^{p}\leq1}%
{\mathbb{E}}_{P}\left[  \int_{0}^{1}\nabla{\mathfrak{h}}(X+t\delta^{1/p}%
\bar{\Delta})\cdot\bar{\Delta}dt]\geq{\mathbb{E}}_{P}^{\left(  p-1\right)
/p}[\left\Vert \nabla{\mathfrak{h}}(X)\right\Vert _{\ast}^{p/\left(
p-1\right)  }\right]  +o\left(  1\right)
\]
uniformly over $P\in{\mathfrak{P}}(\mathcal{[}-\nu,\nu\mathcal{]}^{d})$ and
${\mathbb{E}}_{P}\left\Vert \nabla{\mathfrak{h}}\left(  Y\right)  \right\Vert
\geq b_{0}$. For the upper bound, we can apply the duality representation,
just as we did in the phi-divergence case. Using duality, we have that
\begin{align*}
&  \sup_{{\mathbb{E}}_{P}\left\Vert \bar{\Delta}\right\Vert ^{p}\leq
1}{\mathbb{E}}_{P}\left[  \int_{0}^{1}\nabla{\mathfrak{h}}(X+t\delta^{1/p}%
\bar{\Delta})\cdot\bar{\Delta}dt\right] \\
&  =\min_{\bar{\lambda}>0}\left\{  \bar{\lambda}+{\mathbb{E}}_{P}\left[
\sup_{\bar{\Delta}}\int_{0}^{1}\nabla{\mathfrak{h}}(X+t\delta^{1/p}\bar
{\Delta})\cdot\bar{\Delta}dt-\bar{\lambda}\left\Vert \bar{\Delta}\right\Vert
^{p}\right]  \right\}  .
\end{align*}
Once again, we select a specific choice $\bar{\lambda}_{0}$ given by
\[
0<\bar{\lambda}_{0}=\arg\min\left\{  \bar{\lambda}+{\mathbb{E}}_{P}[\sup
_{\bar{\Delta}}\{\left\Vert \nabla{\mathfrak{h}}(X)\right\Vert _{\ast}%
\cdot\left\Vert \bar{\Delta}\right\Vert -\bar{\lambda}\left\Vert \bar{\Delta
}\right\Vert ^{p}\}]:\bar{\lambda}\geq0\right\}  .
\]
The fact that $\bar{\lambda}_{0}>0$ follows because ${\mathbb{E}}%
_{P}\left\Vert \nabla{\mathfrak{h}}(X)\right\Vert _{\ast}>0$. We then obtain
\begin{align*}
&  \sup_{{\mathbb{E}}_{P}\left\Vert \bar{\Delta}\right\Vert ^{p}\leq
1}{\mathbb{E}}_{P}\left[  \int_{0}^{1}\nabla{\mathfrak{h}}(X+t\delta^{1/p}%
\bar{\Delta})\cdot\bar{\Delta}\right] \\
&  \leq\bar{\lambda}_{0}+{\mathbb{E}}_{P}\left[  \sup_{\bar{\Delta}}\{\int%
_{0}^{1}\nabla{\mathfrak{h}}(X+t\delta^{1/p}\bar{\Delta})\cdot\bar{\Delta
}dt-\bar{\lambda}_{0}\left\Vert \bar{\Delta}\right\Vert ^{p}\}\right]  .
\end{align*}
Next, we argue that the family of functions%
\[
s_{\delta}\left(  x\right)  :=\sup_{\bar{\Delta}}\left[  \int_{0}^{1}%
\nabla{\mathfrak{h}}(x+t\delta^{1/p}\bar{\Delta})\cdot\bar{\Delta}%
dt-\bar{\lambda}_{0}\left\Vert \bar{\Delta}\right\Vert ^{p}\right]
\]
converges uniformly on compact sets to the function $s_{0}\left(  x\right)  $.
Let us consider the sup over $\left\Vert \bar{\Delta}\right\Vert
>\varepsilon/\delta^{1/p}$ and note because (\ref{Loc_bdd}) is bounded on
compact sets, there exists a constant $c_{0}$ independent of $x\in
\lbrack-\nu,\nu]^{d}$ such that
\begin{align*}
&  \int_{0}^{1}\left(  \nabla{\mathfrak{h}}(x+t\delta^{1/p}\bar{\Delta
})-\nabla{\mathfrak{h}}(x)\right)  \cdot\bar{\Delta}dt-\bar{\lambda}%
_{0}\left\Vert \bar{\Delta}\right\Vert ^{p}\\
&  \leq c_{0}(1+\delta^{\left(  p-1\right)  /p}\left\Vert \bar{\Delta
}\right\Vert ^{p-1})\left\Vert \Delta\right\Vert -\bar{\lambda}_{0}\left\Vert
\bar{\Delta}\right\Vert ^{p}.
\end{align*}
By selecting $\delta$ small enough (depending only on $c_{0}>0$ and
$\bar{\lambda}_{0}>0$) we see that the right hand side can be made arbitrarily
negative uniformly over $x\in\lbrack-\nu,\nu]^{d}$ as $\delta\rightarrow0$. So, it
suffices to consider only $\left\Vert \bar{\Delta}\right\Vert \leq
\varepsilon/\delta^{1/p}$, in this case, since $\nabla{\mathfrak{h}}(\cdot)$
is continuous, then it is uniformly continuous on compacts. So, we can write,
in terms of the (uniform) modulus of the continuity function ${\mathfrak{m}%
}\left(  \cdot\right)  $
\[
\left\Vert \nabla{\mathfrak{h}}(x+t\delta^{1/p}\bar{\Delta})-\nabla
{\mathfrak{h}}(x)\right\Vert \leq{\mathfrak{m}}\left(  \varepsilon\right)  ,
\]
where ${\mathfrak{m}}\left(  \varepsilon\right)  \rightarrow0$. In conclusion,
we have that
\begin{align*}
&  \sup_{\left\Vert \Delta\right\Vert \leq\varepsilon/\delta^{1/p}}%
[\nabla{\mathfrak{h}}(x)\cdot\bar{\Delta}\left(  1-{\mathfrak{m}}\left(
\varepsilon\right)  \right)  -\bar{\lambda}_{0}\left\Vert \bar{\Delta
}\right\Vert ^{p}]\\
&  \leq s_{\delta}\left(  x\right)  \leq\sup_{\left\Vert \Delta\right\Vert
\leq\varepsilon/\delta^{1/p}}[\nabla{\mathfrak{h}}(x)\cdot\bar{\Delta}\left(
1+{\mathfrak{m}}\left(  \varepsilon\right)  \right)  -\bar{\lambda}%
_{0}\left\Vert \bar{\Delta}\right\Vert ^{p}].
\end{align*}
Further, the range $\left\Vert \Delta\right\Vert \leq\varepsilon/\delta^{1/p}$
in the upper and lower envelopes above can be further constrained to be
compact (independent of $\varepsilon$ and $\delta$, but depending on
$\bar{\lambda}_{0}>0$). From the above expressions, we deduce the required
uniform convergence of $s_{\delta}\left(  \cdot\right)  \rightarrow
s_{0}\left(  \cdot\right)  $ on compacts. The asymptotic upper bound then
follows from these estimates.\textrm{$\hfill\square$} \newline

Similar results have appeared in the literature (cf., \cite{BDOW21}). An important
difference which is useful in our analysis is that the above result is uniform
over a class $P\in{\mathfrak{P}}(\mathcal{[}-\nu,\nu\mathcal{]}^{d})$ such that 
${\mathbb{E}}_{P}\left\Vert \nabla{\mathfrak{h}}\left(  Y\right)  \right\Vert
\geq b_{0}$.

In order to write the expansion of $\mathcal{F}_{n}(\theta,\delta_{n})$ we
clarify that here we use $\nabla_{x}l\left(  x,\theta\right)  $ to denote the
gradient with respect to $x$. Under suitable boundedness and smoothness
assumptions, the previous result yields

\begin{equation}
\mathcal{F}_{n}(\theta,\delta_{n})=f_{n}\left(  \theta\right)  +\delta
_{n}^{1/p}{\mathbb{E}}_{P_{n}}^{\left(  p-1\right)  /p}[\left\Vert \nabla
_{x}l\left(  X,\theta\right)  \right\Vert _{\ast}^{p/\left(  p-1\right)
}]+\delta_{n}^{1/p}\epsilon_{n}\left(  \theta\right)  . \label{Wass_F_VAl}%
\end{equation}
We collect the precise statement of our result next. The proof is similar to
that of Theorem \ref{th-asymdro} and thus omitted.

\begin{theorem}
\label{th_value_DRO_W}Suppose $l\left(  \cdot,\theta\right)  $ is continuously
differentiable, that
\begin{equation}
\left(  x,\theta\right)  \mapsto\sup\{\left\Vert \nabla l\left(
x+\Delta,\theta\right)  -\nabla l\left(  x,\theta\right)  \right\Vert
/(1+\left\Vert \Delta\right\Vert ^{p-1}):\left\Vert \Delta\right\Vert \geq0\}
\label{eLB}%
\end{equation}
is locally bounded, that $P_{\ast}$ has compact support, $l\left(
x,\cdot\right)  $ is Lipschitz continuous and
\begin{equation}
\inf_{\theta\in\Theta^{\ast}}{\mathbb{E}}_{P_{\ast}}[\left\Vert \nabla
_{x}l\left(  X,\theta\right)  \right\Vert ]>0. \label{eNon_Deg}%
\end{equation}
Then, we have the following types of behavior of optimal values .\newline%
 {\rm (a-W)} If $\delta_{n}^{1/p}=o\left(  n^{-1/2}\right)  $, then%
\begin{equation}
n^{1/2}\left(  \bar{\vartheta}_{n}-\vartheta\right)  \dst\min_{\theta\in
\Theta^{\ast}}{\mathfrak{g}}(\theta).
\end{equation}
 {\rm (b-W)} If $\delta_{n}^{1/p}=\beta n^{-1/2}$, then
\[
n^{1/2}\left(  \bar{\vartheta}_{n}-\vartheta\right)  \dst\min_{\theta\in
\Theta^{\ast}}\left\{  {\mathfrak{g}}(\theta)+{\mathbb{E}}_{P_{\ast}}^{\left(
p-1\right)  /p}\big[\left\Vert \nabla_{x}l\left(  X,\theta\right)  \right\Vert
_{\ast}^{p/\left(  p-1\right)  }\big ]\right\}  .
\]
 {\rm (c-W)} If $o\left(  \delta_{n}^{1/p}\right)  =n^{-1/2}$, then
\[
\delta_{n}^{-1/p}\left(  \bar{\vartheta}_{n}-\vartheta\right)  \dst\min
_{\theta\in\Theta^{\ast}}{\mathbb{E}}_{P_{\ast}}^{\left(  p-1\right)
/p}[\left\Vert \nabla_{x}l\left(  X,\theta\right)  \right\Vert _{\ast
}^{p/\left(  p-1\right)  }].
\]

\end{theorem}

A completely analogous situation to Example \ref{ex-abs} can also be
constructed in Wasserstein DRO to show that both differentiability of
$l\left(  \cdot,\theta\right)  $ and $p>1$ are important in deriving the
asymptotics in the critical case. The case in which $p=2$ was covered in
\cite{BMS22} under suitable quadratic growth conditions and the existence of
second moments.

\setcounter{equation}{0}

\section{General Principle in Action: Optimal Solutions}

\label{Section_GP_OS}

We complete our discussion in this section, considering optimal solutions for
phi-divergence and DRO problems. A key observation is that in both the
phi-divergence case and the Wasserstein DRO case the uncertainty set is
compact in the weak topology and therefore, if Assumption \ref{ass-2} holds, the  function $\mathcal{F}_{n}(\cdot,\delta_{n})$ is differentiable and its gradient    has expansion \eqref{DROGCLT}. In fact, the derivative can be shown to exist if we are able
to argue that, for $\delta$ sufficiently small, the worst case measure is
unique. This is precisely the strategy that we will pursue in this section.
Throughout the section we impose the condition that $\Theta^{\ast}=\left\{
\theta^{\ast}\right\}  $. Recall that $\sigma^{2}(\theta):=\mathrm{Var}%
_{P_{\ast}}(l(X,\theta))$.

\subsection{The Phi-Divergence Case}

\begin{theorem}
\label{thm_DRO_OS_phi} Suppose that Assumptions {\rm \ref{ass-sec}}, {\rm \ref{ass-2}}  and {\rm \ref{ass-phi}} hold, that $l\left(  x,\cdot\right)  $ is essentially bounded under $P_{\ast
}$ and $\sigma^{2}(\theta^{\ast})>0$, and that $\theta_{\ast}(v)$ is directionally differentiable at $v=0$ (in the Hadamard
sense).
Let $Z\sim N(0,\Sigma)$, where $\Sigma$ is the covariance matrix of
$\nabla l(X,\theta^{\ast})$. 
Then we have  the following.
\newline {\rm (A-phi)} If $\delta
_{n}=o\left(  n^{-1}\right)  $, then%
\[
n^{1/2}\left(  \bar{\theta}_{n}-\theta_{\ast}\right)  \dst\theta_{\ast
}^{\prime}\left(  0,Z\right)  .
\]
 {\rm (B-phi)} If $\delta_{n}=\beta n^{-1}$, then
\[
n^{1/2}\left(  \bar{\theta}_{n}-\theta_{\ast}\right)  \dst\,\theta_{\ast
}^{\prime}\left(  0,Z+\kappa^{1/2}\beta^{1/2}\nabla\sigma(\theta^{\ast
})\right)  .
\]
 {\rm (C-phi)} If $o\left(  \delta_{n}\right)  =n^{-1}$, then
\[
\delta_{n}^{-1/2}\left(  \bar{\theta}_{n}-\theta_{\ast}\right)  \dst\,\theta
_{\ast}^{\prime}\left(  0,\kappa^{1/2}\beta^{1/2}\nabla\sigma(\theta^{\ast
})\right)  .
\]

\end{theorem}

\textbf{Proof.} Applying the centering and scaling used to obtain
(\ref{OPT_PHI}) we obtain
\[
\mathcal{F}_{n}(\theta,\delta_{n})=f_{n}\left(  \theta\right)  +\delta
_{n}^{1/2}D_{n}\left(  \theta,\delta_{n}\right)  ,
\]
where
\begin{equation}
\mathcal{D}_{n}\left(  \theta,\delta_{n}\right)  =\sup_{{\mathbb{E}}_{P_{n}%
}\left(  \Delta\right)  =0,\Delta\geq-\delta_{n}^{-1/2},\delta_{n}%
^{-1/2}{\mathbb{E}}_{P_{n}}(\phi(1+\delta_{n}^{1/2}\Delta))\leq1}{\mathbb{E}%
}_{P_{n}}[l_{n}\left(  X,\theta\right)  \Delta], \label{eDprimal}%
\end{equation}
and
\[
\bar{l}_{n}\left(  X,\theta\right)  =l\left(  X,\theta\right)  -f_{n}\left(
\theta\right)  .
\]
It suffices to show that
\[
\nabla\mathcal{D}_{n}\left(  \theta,\delta_{n}\right)  \rightarrow
\nabla\varrho\left(  \theta\right)
\]
uniformly over some region $\left\Vert \theta-\theta_{\ast}\right\Vert
\leq\delta_{0}$ for some $\delta_{0}$. Note that the optimization region in
(\ref{eDprimal}) is compact in the weak topology and therefore, by Danskin's
Theorem (see \cite[sections 5.1.3 and 7.1.5]{SDR}, Section 7), we have that
$\mathcal{D}_{n}\left(  \cdot,\delta_{n}\right)  $ is directionally
differentiable and by the uniqueness of the optimal $\bar{\Delta}_{n}$ for
$\delta_{n}$ sufficiently small we have that
\[
\nabla\mathcal{D}_{n}\left(  \theta,\delta_{n}\right)  ={\mathbb{E}}_{P_{n}%
}[\nabla l_{n}\left(  X,\theta\right)  \bar{\Delta}_{n}\left(  \theta\right)
].
\]
We can precisely characterize $\bar{\Delta}_{n}\left(  \theta\right)  $ from
Proposition \ref{Prop_Expand_Phi} over a region $\left\Vert \theta
-\theta_{\ast}\right\Vert \leq\delta_{0}$ for which we can guarantee
$Var_{P_{n}}[l\left(  X,\theta\right)  ]>0$. Note that such $\delta_{0}>0$ can
be found assuming that $n>N$ (for some random but finite almost surely $N $
because of the Strong Law of Large Numbers and continuity since $Var_{P_{\ast
}}[l\left(  X,\theta_{\ast}\right)  ]>0$. We have, uniformly over $\left\Vert
\theta-\theta_{\ast}\right\Vert \leq\delta_{0}$, for $n>N$,%
\[
\bar{\Delta}_{n}\left(  \theta\right)  =\sqrt{\kappa}\frac{l_{n}\left(
X,\theta\right)  }{\sqrt{\phi^{\prime\prime}\left(  1\right)  Var_{P_{n}%
}[l\left(  X,\theta\right)  ]}}+\epsilon_{n}\left(  \theta\right)  .
\]
On the other hand, defining
\[
\bar{l}\left(  X,\theta\right)  =l\left(  X,\theta\right)  -f\left(
\theta\right)  ,
\]
we have that
\[
\nabla\varrho\left(  \theta\right)  ={\mathbb{E}}_{P_{\ast}}[\bar{l}%
(X,\theta)\cdot\bar{\Delta}\left(  \theta\right)  ],
\]
where
\[
\bar{\Delta}\left(  \theta\right)  =\sqrt{\kappa}\frac{\bar{l}\left(
X,\theta\right)  }{\sqrt{\phi^{\prime\prime}\left(  1\right)  Var_{P_{\ast}%
}[l\left(  X,\theta\right)  ]}}.
\]
We obtain
\begin{align*}
&  \nabla\mathcal{D}_{n}\left(  \theta,\delta_{n}\right)  -\nabla
\varrho\left(  \theta\right) \\
&  ={\mathbb{E}}_{P_{n}}[\nabla l_{n}\left(  X,\theta\right)  \bar{\Delta}%
_{n}\left(  \theta\right)  ]-{\mathbb{E}}_{P_{\ast}}[\nabla\bar{l}%
(X,\theta)\cdot\bar{\Delta}\left(  \theta\right)  ]\\
&  ={\mathbb{E}}_{P_{n}}[\left(  \nabla l_{n}\left(  X,\theta\right)
-\nabla\bar{l}(X,\theta)\right)  \bar{\Delta}_{n}\left(  \theta\right)  ]\\
&  +{\mathbb{E}}_{P_{n}}[\nabla\bar{l}(X,\theta)(\bar{\Delta}_{n}\left(
\theta\right)  -\bar{\Delta}\left(  \theta\right)  )]\\
&  +{\mathbb{E}}_{P_{n}}[\nabla\bar{l}\left(  X,\theta\right)  \bar{\Delta
}\left(  \theta\right)  ]-{\mathbb{E}}_{P_{\ast}}[\nabla\bar{l}(X,\theta
)\cdot\bar{\Delta}\left(  \theta\right)  ].
\end{align*}
It follows that
\[
\bar{\Delta}_{n}\left(  \theta\right)  \rightarrow\bar{\Delta}\left(
\theta\right)
\]
uniformly over $\left\Vert \theta-\theta_{\ast}\right\Vert \leq\delta_{0}$,
and
\[
\left(  \nabla l_{n}\left(  X,\theta\right)  -\nabla\bar{l}(X,\theta)\right)
\rightarrow0
\]
uniformly in probability (in fact almost surely) as $n\rightarrow0$. Uniform
convergence in probability over $\left\Vert \theta-\theta_{\ast}\right\Vert
\leq\delta_{0}$ follows from these observations.\textrm{$\hfill\square$}

\subsection{The Wasserstein Distance Case}

As in Proposition \ref{Prof_Sens_WDRO}, in order to simplify the exposition,
we assume that $P_{\ast}$ has a compact support. We also let $\left\Vert
\cdot\right\Vert $ be the $\ell_{\bar{p}}$ norm for $\bar{p}\in\left(
1,\infty\right)  $. This choice, in particular, satisfies that for any $x$
such that $\left\Vert x\right\Vert =1$, the set%
\begin{equation}
\arg\max\{z^{T}x:\left\Vert z\right\Vert =1\}\text{ is a singleton.}
\label{eq_singleton}%
\end{equation}
This will help us argue, in the presence of Lipschitz gradients, that the
worst case adversarial distribution is unique when the distributional
uncertainty, $\delta$, is sufficiently small and this, in turn, will help
guarantee differentiability. In this section, we use $\nabla_{\theta}$ and
$\nabla_{x}$ to denote the derivatives with respect to $\theta$ and $x$,
respectively. The derivative with respect to all of the arguments is simply
denoted via $\nabla$.

\begin{theorem}
\label{Thm_DRO_OS_W} Suppose that Assumptions {\rm \ref{ass-sec}} and {\rm \ref{ass-2}} hold.
Further, assume that
conditions \eqref{eLB} - \eqref{eNon_Deg} hold and that $\theta_{\ast}(v)$ is
directionally differentiable at $v=0$ (in the Hadamard sense). Let $Z\sim
N(0,\Sigma)$, where $\Sigma$ is the covariance matrix of $\nabla
l(X,\theta^{\ast})$. Then, we have the following.\newline
 {\rm (A-W)} If
$\delta_{n}^{1/p}=o\left(  n^{-1/2}\right)  $, then%
\begin{equation}
n^{1/2}\left(  \bar{\theta}_{n}-\theta_{\ast}\right)  \dst\theta_{\ast
}^{\prime}\left(  0,Z\right)  .
\end{equation}
 {\rm (B-W)} If $\delta_{n}^{1/p}=\beta n^{-1/2}$, then
\begin{equation}
n^{1/2}\left(  \bar{\theta}_{n}-\theta_{\ast}\right)  \dst\theta_{\ast
}^{\prime}\left(  0,Z+\nabla_{\theta}{\mathbb{E}}_{P_{\ast}}^{\left(
p-1\right)  /p}[\left\Vert \nabla_{x}l\left(  X,\theta_{\ast}\right)
\right\Vert _{\ast}^{p/\left(  p-1\right)  }]\right)  .
\end{equation}
 {\rm (C-W)} If $o\left(  \delta_{n}^{1/p}\right)  =n^{-1/2}$, then
\begin{equation}
\delta_{n}^{-1/2}\left(  \bar{\theta}_{n}-\theta_{\ast}\right)  \dst\theta
_{\ast}^{\prime}\left(  0,\nabla_{\theta}{\mathbb{E}}_{P_{\ast}}^{\left(
p-1\right)  /p}[\left\Vert \nabla_{x}l\left(  X,\theta_{\ast}\right)
\right\Vert _{\ast}^{p/\left(  p-1\right)  }]\right)  .
\end{equation}

\end{theorem}

\textbf{Proof.} Most of the work has already been done in the proof of
Proposition \ref{Prof_Sens_WDRO}. We have that
\[
\mathcal{F}_{n}(\theta,\delta_{n})=f_{n}\left(  \theta\right)  +\delta
_{n}^{1/p}\mathcal{D}_{n}\left(  \theta,\delta_{n}\right)  ,
\]
where
\begin{equation}
\mathcal{D}_{n}\left(  \theta,\delta_{n}\right)  :=\sup_{{\mathbb{E}}_{P_{n}%
}\left\Vert \bar{\Delta}\right\Vert ^{p}\leq1}{\mathbb{E}}_{P_{n}}\left[
\int_{0}^{1}\nabla_{x}l(X+t\delta_{n}^{1/p}\bar{\Delta},\theta)\cdot
\bar{\Delta}dt\right]  . \label{Good_for_Diff}%
\end{equation}
It suffices to show uniform convergence of $\nabla\mathcal{D}_{n}\left(
\theta,\delta_{n}\right)  $ to $\nabla\varrho\left(  \theta\right)  $ in some
neighborhood $\left\Vert \theta-\theta_{\ast}\right\Vert \leq\delta_{0}$ for
some $\delta_{0}>0$.

From the proof of Proposition \ref{Prof_Sens_WDRO} we can collect several
facts, note that we are assuming that $\nabla l(\cdot)$ is $L$-Lipschitz,
which guarantees (\ref{eLB}).

I) First,
\[
\sup_{\theta\in\Theta}\left\vert \mathcal{D}_{n}\left(  \theta,\delta
_{n}\right)  -\varrho\left(  \theta\right)  \right\vert \rightarrow0
\]
in probability.

II)\ Moreover, we also saw that there exists a random $N$ (finite with
probability one) such that
\begin{align*}
&  \mathcal{D}_{n}\left(  \theta,\delta_{n}\right) \\
&  =\bar{\lambda}_{n}\left(  \theta\right)  +{\mathbb{E}}_{P_{n}}\left[
\max_{\bar{\Delta}}\int_{0}^{1}\nabla_{x}l(X+t\delta_{n}^{1/p}\bar{\Delta
},\theta)\cdot\bar{\Delta}dt-\bar{\lambda}_{n}\left(  \theta\right)
\left\Vert \bar{\Delta}\right\Vert ^{p}\right]  ,
\end{align*}
with $\bar{\lambda}_{n}\left(  \theta\right)  >0$ for all $n>N$ uniformly over
$\left\Vert \theta-\theta_{\ast}\right\Vert \leq\delta_{0}$ for $\delta_{0}>0$
small enough so that ${\mathbb{E}}_{P_{\ast}}\left[  \left\Vert \nabla
l(X,\theta)\right\Vert \right]  >\delta_{0}$. Note that such $\delta_{0}$
exists by continuity since we assume that ${\mathbb{E}}_{P_{\ast}}\left[
\left\Vert \nabla l(X,\theta_{\ast})\right\Vert \right]  >0$.

III)\ Finally, also on the set $\left\Vert \theta-\theta_{\ast}\right\Vert
\leq\delta_{0}$ from II), since $\nabla l(\cdot)$ is Lipschitz for all
$\delta_{n}$ sufficiently small, we have that the maximizer $\bar{\Delta}%
_{n}\left(  X,\theta\right)  $ inside the expectation is unique because of
(\ref{eq_singleton}) and it converges uniformly on compacts in the variable
$X$ and over $\left\Vert \theta-\theta_{\ast}\right\Vert \leq\delta_{0}$ to
\begin{equation}
\bar{\Delta}\left(  X,\theta\right)  =\bar{\Delta}_{1}\left(  X,\theta\right)
\left\Vert \nabla l(X,\theta)\right\Vert _{\ast}^{1/\left(  p-1\right)
}{\mathbb{E}}_{P}^{-1/p}\left\Vert \nabla l(X,\theta)\right\Vert _{\ast
}^{p/\left(  p-1\right)  }, \label{Def_Del}%
\end{equation}
where%

\[
\bar{\Delta}_{1}\left(  X,\theta\right)  =\arg\max\{\nabla l(X,\theta
)\cdot\bar{\Delta}:\left\Vert \bar{\Delta}\right\Vert =1\}.
\]

Next, by Danskin's Theorem, see \cite[sections 5.1.3 and 7.1.5]{SDR}, Section
7, because the uncertainty set is compact in the weak topology, we have that
$\mathcal{D}_{n}\left(  \cdot,\delta_{n}\right)  $ is differentiable by
uniqueness of $\bar{\Delta}_{n}\left(  X,\theta\right)  $. The most convenient
representation to see this is (\ref{Good_for_Diff}). It is also direct that
$\varrho\left(  \cdot\right)  $ is differentiable everywhere. Moreover, since
\[
\varrho\left(  \theta\right)  =\sup_{{\mathbb{E}}_{P_{\ast}}\left\Vert
\Delta\right\Vert ^{p}\leq1}{\mathbb{E}}_{P_{\ast}}[\nabla_{x}l(X,\theta
)\cdot\Delta],
\]
Danskin's Theorem also applies and we have that
\[
\nabla_{\theta}\varrho\left(  \theta\right)  ={\mathbb{E}}_{P_{\ast}}%
[\nabla_{\theta,x}l(X,\theta)\cdot\bar{\Delta}],
\]
where $\bar{\Delta}$ is given in (\ref{Def_Del}). So, we have (using
$\bar{\Delta}_{n}$) instead of $\bar{\Delta}_{n}\left(  X,\theta\right)  $,
\begin{align*}
&  \nabla_{\theta}\mathcal{D}_{n}\left(  \theta,\delta_{n}\right)
-\nabla_{\theta}\varrho\left(  \theta\right) \\
&  ={\mathbb{E}}_{P_{n}}\left[  \int_{0}^{1}\nabla_{\theta,x}l(X+t\delta
_{n}^{1/p}\bar{\Delta}_{n},\theta)\cdot\bar{\Delta}_{n}dt]-{\mathbb{E}%
}_{P_{\ast}}[\nabla_{\theta,x}l(X,\theta)\cdot\bar{\Delta}\right] \\
&  ={\mathbb{E}}_{P_{n}}\left[  \int_{0}^{1}\left(  \nabla_{\theta
,x}l(X+t\delta_{n}^{1/p}\bar{\Delta}_{n},\theta)-\nabla_{\theta,x}%
l(X,\theta)\right)  \cdot\bar{\Delta}_{n}dt\right] \\
&  +{\mathbb{E}}_{P_{n}}[\nabla_{\theta,x}l(X,\theta)\cdot\bar{\Delta}%
_{n}]-{\mathbb{E}}_{P_{\ast}}[\nabla_{\theta,x}l(X,\theta)\cdot\bar{\Delta}].
\end{align*}
Since $\bar{\Delta}_{n}\left(  X,\cdot\right)  $ converges (in probability)
uniformly over compact sets in $X$ and $\left\Vert \theta-\theta_{\ast
}\right\Vert \leq\delta_{0}$ and it is bounded almost surely, we obtain the
required uniform convergence occurs in probability from the fact that $\nabla
l(\cdot)$ is Lipschitz continuous.\textrm{$\hfill\square$} \newline

A similar result is obtained in \cite{BMS22} quadratic growth conditions and
the existence of second moments (thus relaxing compactness assumptions).
However, \cite{BMS22} primarily focuses on the case in which the optimal
solution lies in the interior of the feasible region. Our discussion here can
be used in combination with the analysis in \cite{BMS22} to deal with boundary cases.
\\
\\
\textbf{Acknowledgement}\\
 J. Blanchet's research was partially
supported by the Air Force Office of Scientific Research (AFOSR) under award
number FA9550-20-1-0397, with additional support from NSF 1915967 and 2118199.
The research of A. Shapiro was partially supported by (AFOSR) Grant FA9550-22-1-0244.

\bibliographystyle{plain}
\bibliography{references,references-Jose}

\end{document}